\documentclass[11pt]{amsart}
\usepackage{amssymb}

\newtheorem{definition}{Definition}[section]
\newtheorem{theorem}[definition]{Theorem}
\newtheorem{lemma}[definition]{Lemma}
\newtheorem{proposition}[definition]{Proposition}

\begin{document}

\title{The hyperoctahedral quantum group}
\author{Teodor Banica}
\address{T.B.: Department of Mathematics,
Paul Sabatier University, 118 route de Narbonne, 31062 Toulouse,
France} \email{banica@picard.ups-tlse.fr}
\author{Julien Bichon}
\address{J.B.: Department of Mathematics, Blaise Pascal University, Campus des Cezeaux, 63177 Aubiere Cedex, France}
\email{Julien.Bichon@math.univ-bpclermont.fr}
\author{Beno\^it Collins}
\address{B.C.: Department of Mathematics,
Claude Bernard University, 43 bd du 11 novembre 1918, 69622
Villeurbanne, France, and}
\address{\vskip -6mm Department of Mathematics, University of Ottawa, 585 King Edward, Ottawa, ON K1N 6N5, Canada} \email{collins@math.univ-lyon1.fr}

\begin{abstract}
We consider the hypercube in $\mathbb R^n$, and show that its quantum symmetry group is a $q$-deformation of $O_n$ at $q=-1$. Then we consider the graph formed by $n$ segments, and show that its quantum symmetry group is free in some natural sense. This latter quantum group, denoted $H_n^+$, enlarges Wang's series $S_n^+,O_n^+,U_n^+$.
\end{abstract}

\maketitle

\section*{Introduction}

The idea of noncommuting coordinates goes back to Heisenberg. Several theories emerged
from Heisenberg's work, most complete being Connes' noncommutative
geometry, where the base space is a Riemannian manifold. See \cite{cn}.

A natural question is about studying algebras of free coordinates
on algebraic groups. Given a group $G\subset U_n$, the matrix
coordinates $u_{ij}\in C(G)$ commute with each other, and satisfy
certain relations $R$. One can define then the universal algebra
generated by abstract variables $u_{ij}$, subject to the relations
$R$. The spectrum of this algebra is an abstract object, called
noncommutative version of $G$. The noncommutative version is not
unique, because it depends on $R$. We have the following examples:
\begin{enumerate}
\item Free quantum semigroups. The algebra generated by variables
$u_{ij}$ with relations making $u=(u_{ij})$ a unitary matrix was
considered by Brown \cite{br}. This algebra has a comultiplication
and a counit, but no antipode. In other words, the corresponding
noncommutative version $U_n^{nc}$ is a quantum semigroup. \item
Quantum groups. A remarkable discovery is that for a compact Lie
group $G$, once commutativity is included into the relations $R$,
one can introduce a complex parameter $q$, as to get deformed
relations $R^q$. The corresponding noncommutative versions $G^q$
are called quantum groups. See Drinfeld \cite{dr}. \item Compact
quantum groups. The quantum group $G^q$ is semisimple for $q$ not
a root of unity, and compact under the slightly more restrictive
assumption $q\in\mathbb R$. Woronowicz gave in \cite{wo1} a simple
list of axioms for compact quantum groups, which allows
construction of many other noncommutative versions.\item Free
quantum groups. The Brown algebras don't fit into
Woronowicz's axioms, but a slight modification leads to free
quantum groups. The quantum groups $O_n^+,U_n^+$ appeared in
Wang's thesis \cite{wa1}. Then Connes suggested use of symmetric
groups, and the quantum group $S_n^+$ was constructed in
\cite{wa2}.
\end{enumerate}

The quantum groups $S_n^+,O_n^+,U_n^+$ have been intensively studied in the last years.

The purpose of this paper is to introduce and study a new free quantum group. This is $H_n^+$, the free analogue of the hyperoctahedral group $H_n$.

The various algebraic results in \cite{ba2}, \cite{bb1}, \cite{bi2} and integration techniques in \cite{bc1}, \cite{bc2}, \cite{bc3} apply to $H_n^+$, and our results can be summarized in the following table:

\begin{center}
\begin{tabular}[t]{|l|l|l|l|l|}
\hline &$S_n^+$&$H_n^+$&$O_n^+$&$U_n^+$\\
\hline Jones diagrams&$TL$&colored $TL$&double $TL$&free $TL$\\
\hline Speicher partitions&$NC$&colored $NC$&double $NC$&free $NC$\\
\hline Voiculescu law&free Poisson&free Bessel&semicircular&circular\\
\hline Di Francesco matrix&full meander&colored meander&meander&free meander\\
\hline
\end{tabular}
\end{center}
\bigskip

Here the $S_n^+,O_n^+,U_n^+$ columns are from \cite{bc1}, \cite{bc2}, and the $H_n^+$ column is new. The precise meaning of this table, and the connection with \cite{df}, \cite{jo}, \cite{sp}, \cite{vo} and with the paper of Bisch and Jones \cite{bj} will be explained in the body of the paper. Before proceeding, we would like to make a few comments:
\begin{enumerate}
\item The quantum group $H_n^+$ is by definition a subgroup of
$S_{2n}^+$, and emerges from the recent litterature on quantum
permutation groups, via a combination of the approach in
\cite{bb1}, \cite{bb2}, \cite{bbc} with the one in \cite{bc1},
\cite{bc2}, \cite{bc3}.

\item The origins
of $H_n^+$ go back to results in \cite{bi2}, with the construction
of free wreath products, and with the exceptional isomorphism
$H_2^+=O_2^{-1}$. Besides the above table, we have the quite
surprising result that $O_n^{-1}$ is a noncommutative analogue of
$H_n$, for any $n$.

\item The fact that the Fuss-Catalan algebra can be realized with compact quantum groups is known since \cite{ba2}. However, there are many choices for such a quantum group. The results in this paper show that $H_n^+$ is the good choice.

\item The computations for $H_n^+$ belong to a
certain order 1 problematics, for the Fuss-Catalan algebra. For general free wreath products the
situation is discussed in \cite{bb1}, with a conjectural order $0$
statement. The plan here would be to prove and improve this
statement, then to compare it with results of \'Sniady in \cite{sn}.

\item The quantum groups considered in this paper
satisfy $S_n^+\subset G\subset U_n^+$. This condition is very
close to freeness, and we think that quantum groups satisfying it
should deserve a systematic investigation. Some general results in
this sense, at level of integration formulae, are worked out here.
\end{enumerate}

Summarizing, the quantum group $H_n^+$ appears to be a central object of the theory, of the same type as Wang's fundamental examples $S_n^+,O_n^+,U_n^+$.

Finally, let us mention that the semigroup of measures for $H_n^+$, that we call ``free Bessel'', seems to be new. Here are the even moments of these measures:
$$\int x^{2k}\,d\mu_t(x)=\sum_{b=1}^k\frac{1}{b}\begin{pmatrix}k-1\cr b-1\end{pmatrix}\begin{pmatrix}2k\cr b-1\end{pmatrix}t^b$$

As for the odd moments, these are zero. The relation with Bessel functions is that in the classical case, namely for the group $H_n$, the measures that we get are supported on $\mathbb Z$, with density given by Bessel functions of the first kind.

The paper is organized in four parts, as follows:
\begin{enumerate}
\item In 1-2 we discuss the integration problem, with a
presentation of some previously known results, and with an explicit computation for $H_n$. 

\item In 3-4 we show that $O_n^{-1}$ is the quantum symmetry group of the hypercube in $\mathbb R^n$. Thus $O_n^{-1}$ can be regarded as a noncommutative version of $H_n$.

\item In 5-10 we construct the free version $H_n^+$, as quantum symmetry group of the graph formed by $n$ segments. We
perform a detailed combinatorial study of $H_n^+$, the main result being the asymptotic freeness of diagonal coefficients.

\item In 11-12 we extend a part of the
results to a certain class of quantum groups, that we call free.
We give some other examples of free quantum groups, by using a
free product construction, and we discuss the classification
problem.
\end{enumerate}

\subsection*{Acknowledgements} We would like to express our gratitude to the
universities of Lyon, Pau, Toulouse and to the research centers in
Bedlewo and Luminy, for their warm hospitality and support, at
various stages of this project.

\section{The integration problem}

In this section we present the basic idea of Weingarten's
approach to the computation of integrals over compact
groups, and its extension to quantum groups. 

Consider a compact group $G\subset O_n$. This can be either a Lie group, like $O_n$ itself, or just a finite group, like $S_n$, or one of its subgroups.

\begin{definition}
For a compact group $G\subset O_n$, the formula
$$g=\begin{pmatrix}u_{11}(g)&&u_{1n}(g)\cr&\ddots&\cr u_{n1}(g)&&u_{nn}(g)\end{pmatrix}$$
defines $n^2$ functions $u_{ij}:G\to\mathbb R$, that we call
matrix coordinates of $G$.
\end{definition}

These functions are coordinates in the usual sense, when $G$
viewed as a real algebraic variety. We regard them as complex functions, $u_{ij}:G\to\mathbb C$.

Consider now the Haar measure of $G$. The integration problem is
to compute the integrals of various products of matrix
coordinates:
$$I(i,j)=\int_Gu_{i_1j_1}(g)\ldots u_{i_kj_k}(g)\,dg$$

It is convenient at this point to switch to operator algebraic 
language. We regard integration as being an abstract functional on
the algebra $C(G)$:
$$\int :C(G)\to\mathbb C$$

The integration problem is to compute this functional. Indeed, by norm
continuity the functional is uniquely determined on the dense
subalgebra generated by entries of $u$, and by linearity
computation is reduced to that of the above integrals.

\begin{definition}
The integration problem is the computation of
$$I(i,j)=\int u_{i_1j_1}\ldots
u_{i_kj_k}$$ in terms of multi-indices $i=(i_1,\ldots,i_k)$ and
$j=(j_1,\ldots,j_k)$.
\end{definition}

The very first example of such an integral, or rather of a sum of
such integrals, is the following well-known quantity:
$$\int (u_{11}+u_{22}+\ldots +u_{nn})^k=
\# \left\{1\in u^{\otimes k}\right\}$$

Here the number on the right is the multiplicity of $1$ into the
$k$-th tensor power of $u$, whose computation is a fundamental
problem in representation theory. The number on the left is the
integral of the $k$-th power of the character of $u$.

This relation with representation theory provides actually a
solution to the integration problem. The algorithm is as follows:
\begin{enumerate}
\item Decompose the representation $u^{\otimes k}$, as sum of
irreducibles. \item Find a basis $D_k$ for the linear space $F_k$
of fixed vectors of $u^{\otimes k}$.\item Compute the Gram matrix
$G_k$ of this basis $D_k$. \item Compute the Weingarten matrix,
$W_k=G_k^{-1}$. \item Express the orthogonal projection $P_k$ onto
$F_k$ in terms of $W_k$.\item Use the fact that the matrix
coefficients of $P_k$ are the integrals $I(i,j)$.
\end{enumerate}

The fact that this works follows from 5-6. We use here an elementary
result from linear algebra: the projection onto a space having an
orthogonal basis is given by a trivial formula, and in case the
basis is not orthogonal, the Gram matrix (or rather its inverse)
comes into the picture, and leads to a more complicated formula. The last assertion is also a standard fact, coming from Peter-Weyl theory.

The rest of this section contains a detailed
discussion of the problem in 1-4. For the moment, let us summarize the
above discussion, in the form of a vague statement:

\begin{definition}
The Weingarten problem is to find a formula of type
$$\int u_{i_1j_1}\ldots
u_{i_kj_k}=\sum_{pq}\delta_{pi}\delta_{qj}W_{k}(p,q)$$ where the
sum is over diagrams, the $\delta$-symbols encode the coupling of
diagrams with multi-indices, and $W_{k}$ is the inverse of the
coupling matrix of diagrams.
\end{definition}

The origins of this approach go back to Weingarten's paper
\cite{we}. The formula is worked out in detail in \cite{cl},
\cite{cs1}, \cite{cs2}, for several classes of Lie groups.

As already mentioned, the general formula is quite theoretical,
because of serious complexity problems in step 4. In fact, the
Weingarten matrix cannot be computed explicitly in the general
case, and only partial results are available.

Fortunately, we have here a hierarchy of reasonable questions:
\begin{enumerate}
\item Order zero: computation of the law of
$u_{11}+\ldots+u_{nn}$.\item First order: computation of the law
of $u_{11}+\ldots +u_{ss}$, with $1\leq s\leq n$. \item Second
order: integration of carefully chosen quantities, coming from
advanced representation theory, or from theoretical
physics.\end{enumerate}

The list is quite heuristic, and is based on recent work on the
subject, rather than on some precise mathematical statement. The
numbers $0,1,2$ should be thought of as being exponents of
$n^{-1}$ in some related asymptotic expansion. See e.g.
\cite{cgm}.

In this paper we restrict attention to first order problems. It is
convenient to consider arbitrary sums of diagonal entries of $u$,
and to denote them as follows.

\begin{definition}
To any union of intervals $I\subset [0,1]$ we associate the
element
$$u_I=\sum_{s\in nI}u_{ss}$$
of the algebra $C(G)$, where $u$ is the matrix of coefficients of
$G$.
\end{definition}

The example to be kept in mind is with the interval $I=[0,t]$. We get here the element $u_I=u_{11}+\ldots +u_{ss}$ with $s=[tn]$, corresponding to the first order problem.

We will need a few notions from classical and free probability. It
is convenient to use here the abstract framework of $C^*$-algebras
endowed with positive traces:
$$\int:A\to\mathbb C$$

An element of $A$ will be called random variable. The law of a
self-adjoint random variable $a$ is the real probability measure
$\mu$ given by the following formula:
$$\int a^k=\int_{\mathbb R} x^k\,d\mu(x)$$

Recall that two elements $a\in A$ and $b\in B$ are independent in
$A\otimes B$, and that there is an abstract notion of independence
axiomatizing this situation. See \cite{vdn}.

The total length of a union of intervals $I$ is denoted $l(I)$.

\begin{theorem}
Consider the groups $O_n,S_n$.
\begin{enumerate}
\item For $O_n$, the variable $u_I$ is asymptotically Gaussian of
parameter $l(I)$. \item For $S_n$, the variable $u_I$ is asymptotically Poisson of parameter
$l(I)$.\item The variables $u_I,u_J$ with $I\cap
J=\emptyset$ are asymptotically independent.\end{enumerate}
\end{theorem}

This statement is pointed out in \cite{bc2}, in a slightly different form.

We present now a quantum analogue of this result. We use
the following key definition, emerging from Wang's work
\cite{wa1}, \cite{wa2}, which was in turn based on Woronowicz's
axiomatization of compact quantum groups \cite{wo1}, \cite{wo2}.

\begin{definition}
Let $u\in M_n(A)$ be a square matrix over a $C^*$-algebra.
\begin{enumerate}
\item $u$ is called orthogonal if all elements $u_{ij}$ are self-adjoint, and $u$ is unitary. In other words, we must have $u=\bar{u}$
and $u^t=u^{-1}$. \item $u$ is called magic unitary if all elements $u_{ij}$ are projections, and on each row and each column of $u$ these
projections are orthogonal, and sum up to $1$.
\end{enumerate}
\end{definition}

The matrix coordinates of $O_n$ form an orthogonal matrix. The matrix coordinates of $S_n$ are the characteristic functions  $u_{ij}=\chi\{\sigma\in S_n\mid \sigma(j)=i\}$, which form a magic unitary. These remarks and the Gelfand theorem lead to the following formulae:
\begin{eqnarray*}
C(O_n)&=&C^*_{com}\left( u_{ij}\mid u=n\times n\mbox{ {\rm
orthogonal}}\right)\cr C(S_n)&=&C^*_{com} \left( u_{ij}\mid
u=n\times n\mbox{ {\rm magic unitary}}\right)
\end{eqnarray*}

We can proceed now with liberation. We call comultiplication, counit and antipode any
morphisms of $C^*$-algebras of the following type:
\begin{eqnarray*}
\Delta&:&A\to A\otimes A\cr \varepsilon&:&A\to {\mathbb C}\cr
S&:&A\to A^{op}
\end{eqnarray*}

An orthogonal Hopf algebra is a $C^*$-algebra $A$ given with an
orthogonal matrix $u$, such that the following formulae define a
comultiplication, a counit and an antipode:
\begin{eqnarray*} \Delta(u_{ij})&=&\sum u_{ik}\otimes
u_{kj}\cr \varepsilon(u_{ij})&=&\delta_{ij}\cr S(u_{ij})&=&u_{ji}
\end{eqnarray*}

The basic example is the algebra $C(G)$, where $G\subset O_n$ is a
compact group. Here $u$ is the matrix of coordinates of $G$, and the axioms are clear from definitions.

\begin{theorem}
The universal algebras
\begin{eqnarray*}
A_o(n)&=&C^*\left( u_{ij}\mid u=n\times n\mbox{ {\rm
orthogonal}}\right)\cr A_s(n)&=&C^* \left( u_{ij}\mid u=n\times
n\mbox{ {\rm magic unitary}}\right)
\end{eqnarray*}
are orthogonal Hopf algebras.
\end{theorem}

See Wang \cite{wa1}, \cite{wa2}. As a conclusion, we have the following heuristic equalities, where
plus superscripts denote free versions of the corresponding classical
groups:
\begin{eqnarray*}
A_o(n)&=&C(O_n^+)\cr
A_s(n)&=&C(S_n^+)
\end{eqnarray*}

We discuss now the liberation aspect of the operation $G\to G^+$.
There are several interpretations here. The idea that we follow,
pointed out in \cite{bc2}, is that with $n\to\infty$ the integral
geometry of $G$ is governed by independence, while that of $G^+$, by freeness.

Recall that two elements $a\in A$ and $b\in B$
are free in $A*B$, and that there is an abstract notion of
freeness axiomatizing this situation. See \cite{vo}, \cite{vdn}.

\begin{theorem}
Consider the algebras $A_o(n),A_s(n)$.
\begin{enumerate}
 \item For $A_o(n)$, the variable $u_I$ is asymptotically semicircular of parameter $l(I)$. \item For $A_s(n)$, the variable $u_I$ is asymptotically free Poisson of parameter $l(I)$.
\item The variables $u_I,u_J$ with $I\cap J=\emptyset$ are
asymptotically free.
\end{enumerate}
\end{theorem}

This statement appears in \cite{bc2}, in a slightly different
form. The assertions follow from results in \cite{bc1},
\cite{bc2}. This will be explained in detail later on.

\section{The hyperoctahedral group}

In this section we review a few well-known facts regarding the
symmetry group of the cube. The whole material is here in order to
be carefully studied in next sections: several problems will
appear from the simple facts explained below.

\begin{definition}
$H_n$ is the common symmetry group of the following spaces:
\begin{enumerate}
\item The cube ${\rm K}_n\subset \mathbb R^n$, regarded as a
metric space ($2^n$ points). \item The graph $K_n$ which looks
like a cube ($2^n$ vertices, $n2^{n-1}$ edges). \item The space ${\rm
I}_n\subset\mathbb R^n$ formed by the $\pm 1$ points on each axis
($2n$ points).\item The graph $I_n$ formed by $n$ segments ($2n$
vertices, $n$ edges).
\end{enumerate}
\end{definition}

The equivalences 1 $\Leftrightarrow$ 2 and 3 $\Leftrightarrow$ 4 are consequences of the
following general principle: when drawing a geometric object O,
the symmetry group of the graph $O$ that we obtain is the same as
the symmetry group of O.

As for the equivalence 1 $\Leftrightarrow$ 3, this follows from the fact that a rotation of $\mathbb R^n$ leaves invariant the cube  if and
only if it leaves invariant the set of middles of each face.

For each segment $I^i$ consider the element $\tau_i\in H_n$ which
returns $I^i$, and keeps the other segments fixed. The elements
$\tau_1,\ldots,\tau_n$ have order two and commute with each other,
so they generate the product of $n$ copies of $Z_2$, denoted in
this paper $L_n$:
$$L_n=Z_2\times\ldots\times Z_2$$

Now the symmetric group $S_n$ acts as well on $X_n$ by permuting
the segments, and it is routine to check that we have
$H_n=L_n\rtimes S_n$. But this latter crossed product is by
definition the wreath product $Z_2\wr S_n$.

\begin{proposition}
$H_n$ has the following properties: \begin{enumerate}\item We have
$S_n\subset H_n\subset O_n$. \item We have $H_n=L_n\rtimes S_n$.
\item We have $H_n=Z_2\wr S_n$.\end{enumerate}
\end{proposition}

We will see later on that the free analogue $H_n^+$ is a free
wreath product. In terms of the hierarchy of integration problems
explained in previous section, the available results for usual and
free wreath products are as follows:
\begin{enumerate}
\item Second order results for wreath products are obtained in
\cite{sn}. \item Order zero results for free wreath products are
obtained in \cite{bb1}.
\end{enumerate}

For the purposes of this paper, what we actually need are the
first order results, for both $H_n$ and $H_n^+$. Most convenient
here is to present complete proofs for everything. We begin with a
discussion of the first order problem for $H_n$.

\begin{theorem}
The asymptotic law of $u_{11}+\ldots +u_{ss}$ with $s=[tn]$ is
given by
$$\beta_t=e^{-t}\sum_{k=-\infty}^\infty\delta_k
\sum_{p=0}^\infty \frac{(t/2)^{|k|+2p}}{(|k|+p)!p!}$$ where
$\delta_k$ is the Dirac mass at $k\in\mathbb Z$.
\end{theorem}

\begin{proof}
We regard $H_n$ as being the symmetry group of the graph
$I_n=\{I^1,\ldots ,I^n\}$ formed by $n$ segments. The diagonal
coefficients are given by:
$$u_{ii}(g)=\begin{cases}\ 0\ \mbox{ if $g$ moves $I^i$}\\
+1\mbox{ if $g$ fixes $I^i$}\\
-1\mbox{ if $g$ returns $I^i$}
\end{cases}$$

We denote by $\uparrow g,\downarrow g$ the number of segments
among $\{I^1,\ldots ,I^s\}$ which are fixed, respectively returned
by an element $g\in H_n$. With this notation, we have:
$$u_{11}+\ldots+u_{ss}=\uparrow g-\downarrow g$$

We denote by ${\rm P}_n$ probabilities computed over the group
$H_n$. The density of the law of $u_{11}+\ldots+u_{ss}$ at a point
$k\geq 0$ is given by the following formula:
\begin{eqnarray*}
D(k)&=&{\rm P}_n(\uparrow g-\downarrow g=k)\cr
&=&\sum_{p=0}^\infty {\rm P}_n(\uparrow g=k+p, \downarrow g=p)
\end{eqnarray*}

Assume first that we have $t=1$. We use the fact that the
probability of $\sigma\in S_n$ to have no fixed points is
asymptotically $1/e$. Thus the probability of $\sigma\in S_n$ to
have $m$ fixed points is asymptotically $1/(em!)$. In terms of
probabilities over $H_n$, we get:
\begin{eqnarray*}
\lim_{n\to\infty}D(k)&=&\lim_{n\to\infty}\sum_{p=0}^\infty
(1/2)^{k+2p}\begin{pmatrix}k+2p\cr k+p\end{pmatrix} {\rm
P}_n(\uparrow g+\downarrow g=k+2p)\cr &=&\sum_{p=0}^\infty
(1/2)^{k+2p}\begin{pmatrix}k+2p\cr
k+p\end{pmatrix}\frac{1}{e(k+2p)!}\cr &=&
\frac{1}{e}\sum_{p=0}^\infty \frac{(1/2)^{k+2p}}{(k+p)!p!}
\end{eqnarray*}

The general case $0<t\leq 1$ follows by performing some
modifications in the above computation. The asymptotic density is
computed as follows:
\begin{eqnarray*}
\lim_{n\to\infty}D(k)&=& \lim_{n\to\infty}\sum_{p=0}^\infty
(1/2)^{k+2p}\begin{pmatrix}k+2p\cr k+p\end{pmatrix} {\rm
P}_n(\uparrow g+\downarrow g=k+2p) \cr &=&\sum_{p=0}^\infty
(1/2)^{k+2p}\begin{pmatrix}k+2p\cr
k+p\end{pmatrix}\frac{t^{k+2p}}{e^t(k+2p)!}\cr
&=&e^{-t}\sum_{p=0}^\infty \frac{(t/2)^{k+2p}}{(k+p)!p!}
\end{eqnarray*}

Together with $D(-k)=D(k)$, this gives the formula in the
statement.
\end{proof}

\begin{theorem}
The measures $\beta_t$ form a truncated one-parameter semigroup
with respect to convolution.
\end{theorem}

\begin{proof}
Consider the Bessel function of the first kind:
$$f_k(t)=\sum_{p=0}^\infty \frac{t^{|k|+2p}}{(|k|+p)!p!}$$

The measure in the statement is given by the following formula:
$$\beta_t=e^{-t}\sum_{k=-\infty}^\infty
\delta_k\,f_k(t/2)$$

Thus its Fourier transform is given by:
$$F\beta_t(y)=e^{-t}\sum_{k=-\infty}^\infty
e^{ky}\,f_k(t/2)$$

We compute now the derivative with respect to $t$:
$$F\beta_t(y)'=-F\beta_t(y)+\frac{e^{-t}}{2}
\sum_{k=-\infty}^\infty e^{ky}\,f_k'(t/2)$$

On the other hand, the derivative of $f_k$ with $k\geq 1$ is given
by:
\begin{eqnarray*}
f_k'(t)&=&\sum_{p=0}^\infty \frac{(k+2p)t^{k+2p-1}}{(k+p)!p!}\cr
&=&\sum_{p=0}^\infty \frac{(k+p)t^{k+2p-1}}{(k+p)!p!}
+\sum_{p=0}^\infty\frac{p\,t^{k+2p-1}}{(k+p)!p!}\cr
&=&\sum_{p=0}^\infty \frac{t^{k+2p-1}}{(k+p-1)!p!}
+\sum_{p=1}^\infty\frac{t^{k+2p-1}}{(k+p)!(p-1)!}\cr
&=&\sum_{p=0}^\infty \frac{t^{(k-1)+2p}}{((k-1)+p)!p!}
+\sum_{p=1}^\infty\frac{t^{(k+1)+2(p-1)}}{((k+1)+(p-1))!
(p-1)!}\cr &=&f_{k-1}(t)+f_{k+1}(t)
\end{eqnarray*}

This computation works in fact for any $k$, so we get:
\begin{eqnarray*}
F\beta_t(y)'&=&-F\beta_t(y)+\frac{e^{-t}}{2}
\sum_{k=-\infty}^\infty e^{ky} (f_{k-1}(t/2)+f_{k+1}(t/2))\cr
&=&-F\beta_t(y)+\frac{e^{-t}}{2} \sum_{k=-\infty}^\infty
e^{(k+1)y}f_{k}(t/2)+e^{(k-1)y}f_{k}(t/2)\cr
&=&-F\beta_t(y)+\frac{e^{y}+e^{-y}}{2}\,F\beta_t(y)\cr
&=&\left(\frac{e^{y}+e^{-y}}{2}-1\right)F\beta_t(y)
\end{eqnarray*}

Thus the log of the Fourier transform is linear in $t$, and we get
the assertion.
\end{proof}

\section{Quantum automorphisms}

We present now a first approach to the noncommutative analogue of
$H_n$. The idea is to use the fact that $H_n$ is the symmetry
group of the cube.

We begin with a discussion involving finite graphs and finite
metric spaces:
\begin{enumerate}
\item A finite graph, with vertices labelled $1,\ldots ,n$, is
uniquely determined by its adjacency matrix $d\in M_n(0,1)$, given
by $d_{ij}=1$ if $i,j$ are connected, and $d_{ij}=0$ if not. This
matrix has to be symmetric, and null on the diagonal.\item A
finite metric space, with points labelled $1,\ldots ,n$, is
uniquely determined by its distance matrix $d\in M_n(\mathbb
R_+)$. This matrix has to be symmetric, $0$ on the diagonal and $>0$ outside, and
its entries must satisfy the triangle inequality.
\end{enumerate}

We see that in both cases, the object under consideration is determined by a scalar matrix $d\in M_n$. In both cases this matrix is real, symmetric, and null on the diagonal, but for our purposes we will just regard it as an arbitrary complex matrix.

\begin{definition}
Let $X$ be a finite graph, or finite metric space, with matrix
$d\in M_n$.
\begin{enumerate}
\item The quantum symmetry algebra of $X$ is given by
$A_X=A_s(n)/I$, where $I$ is the ideal generated by the $n^2$
relations coming from the equality $du=ud$. \item A Hopf algebra
coacting on $X$ is a quotient of $A_X$. In other words, it is a
quotient of $A_s(n)$ in which the equality $du=ud$ holds.
\end{enumerate}
\end{definition}

As a first example, when $X$ is a set, or a complete graph, or a
simplex, we have $A_X=A_s(n)$. This is because in all these
examples, the non-diagonal entries of $d$ are all equal. In other
words, we have $d=\lambda(F-1)$, where $\lambda$ is a scalar and
$F$ is the matrix filled with $1$'s. But the magic unitarity
condition gives $uF=Fu$, hence $I=(0)$.

The basic example of coaction is that of the algebra $C(G_X)$,
where $G_X$ is the usual symmetry group of $X$. In fact, this is
the biggest commutative algebra coacting on $X$. Or, in other
words, it is the biggest commutative quotient of $A_X$. See
\cite{bb1}.

Recall now from previous section that the cube space ${\rm K}_n$
and the cube graph $K_n$ have the same symmetry group, because of
a general ``drawing principle''. Our first result is that this
principle holds as well for quantum symmetry groups.

\begin{proposition}
The cube space ${\rm K}_n$ and the cube graph $K_n$ have the same
quantum symmetry algebra.
\end{proposition}

\begin{proof}
The distance matrix of the cube is of the following form:
$$d=d_1+\sqrt{2}\,d_2+\sqrt{3}\,d_3+\ldots+
\sqrt{n}\,d_n$$

Here for $l=1,\ldots,n$ the matrix $d_l$ is obtained from $d$ by
replacing all lengths different from $\sqrt{l}$ by the number $0$,
and all lengths equal to $\sqrt{l}$ by the number $1$. This is the
same as the adjacency matrix of the graph $X_l$, having as
vertices the vertices of the cube, and with edges drawn between
pairs of vertices at distance $\sqrt{l}$.

Observe that we have $X_1=K_n$. Now the powers of $d_1$ are
computed by counting loops on the cube, and we get formulae of the
following type:
\begin{eqnarray*}
d_1^2&=&x_{21}\,1_n+x_{22}\,d_2\cr d_1^3&=&x_{31}\,1_n+x_{32}\,
d_2+ x_{33}\,d_3\cr &\ldots&\cr
d_1^n&=&x_{n1}\,1_n+x_{n2}\,d_2+x_{n3}\,d_3+\ldots+x_{nn}\,d_n
\end{eqnarray*}

Here $x_{ij}$ are some positive integers. Now these formulae,
together with the above one for $d$, show that we have an equality
of algebras $<d>=<d_1>$. In particular $u$ commutes with $d$ if
and only if it commutes with $d_1$, and this gives the result.
\end{proof}

We have to find the universal Hopf algebra coacting on the cube.
For this purpose, we use the cube graph $K_n$, and we regard it as
a Cayley graph. Recall from previous section that the direct
product of $n$ copies of $Z_2$ is denoted
$L_n=<\tau_1,\ldots,\tau_n>$.

\begin{lemma}
$K_n$ is the Cayley graph of the group $L_n=<\tau_1,\ldots ,\tau_n>$. Moreover, the eigenvectors and eigenvalues of $K_n$ are given by
\begin{eqnarray*}
v_{i_1\ldots i_n}&=&\sum_{j_1\ldots j_n} (-1)^{i_1j_1
+\ldots+i_nj_n}\tau_1^{j_1}\ldots\tau_n^{j_n}\cr
\lambda_{i_1\ldots i_n}&=&(-1)^{i_1}+\ldots +(-1)^{i_n}
\end{eqnarray*}
when the vector space spanned by vertices of $K_n$ is identified
with $C^*(L_n)$.
\end{lemma}

\begin{proof}
For the first assertion, consider the Cayley graph of $L_n$. The vertices are the elements of $L_n$, which are products of the following form:
$$g=\tau_1^{i_1}\ldots\tau_n^{i_n}$$

The sequence of $0-1$ exponents defining such an element
determines a point of $\mathbb R^n$, which is a vertex of the
cube. Thus the vertices of the Cayley graph are the vertices of the cube.
Now regarding edges, in the Cayley graph these are drawn between
elements $g,h$ having the property $g=h\tau_i$ for some $i$. In
terms of coordinates, the operation $h\to h\tau_i$ means to switch
the sign of the $i$-th coordinate, and to keep the other
coordinates fixed. In other words, we get in this way the edges of
the cube.

For the second part, recall first that the eigenspaces and
eigenvalues of a graph $X$ are by definition those of its
adjacency matrix $d$, and that the action of $d$ on functions on
vertices is by summing over neighbors:
$$df(p)=\sum_{q\sim p}f(q)$$

Now by identifying vertices with elements of $L_n$, hence
functions on vertices with elements of the group algebra
$C^*(L_n)$, we get the following formula:
$$dv=\tau_1v+\ldots +\tau_nv$$

For the vector $v_{i_1\ldots i_n}$ in the statement, we have the
following computation:
\begin{eqnarray*}
dv_{i_1\ldots i_n}&=&\sum_s\tau_s\sum_{j_1\ldots j_n} (-1)^{i_1j_1
+\ldots+i_nj_n}\tau_1^{j_1}\ldots\tau_n^{j_n}\cr
&=&\sum_s\sum_{j_1\ldots j_n} (-1)^{i_1j_1
+\ldots+i_nj_n}\tau_1^{j_1}\ldots\tau_s^{j_s+1}\ldots\tau_n^{j_n}\cr
&=&\sum_s\sum_{j_1\ldots j_n} (-1)^{i_s}(-1)^{i_1j_1
+\ldots+i_nj_n}\tau_1^{j_1}\ldots\tau_s^{j_s}\ldots\tau_n^{j_n}\cr
&=&\sum_s(-1)^{i_s}\sum_{j_1\ldots j_n} (-1)^{i_s}(-1)^{i_1j_1
+\ldots+i_nj_n}\tau_1^{j_1}\ldots\tau_n^{j_n}\cr
&=&\lambda_{i_1\ldots i_n}v_{i_1\ldots i_n}
\end{eqnarray*}

Together with the fact that the vectors $v_{i_1\ldots i_n}$ are
linearly independent, and that there are $2^n$ of them, this gives
the result.
\end{proof}

The above result shows that the eigenvalues of $K_n$ appear as
sums of $\pm 1$ terms, with $n$ terms in each sum. Thus the
eigenvalues are the numbers $s\in\{-n,\ldots ,n\}$. The
corresponding eigenspaces are denoted $E_s$.

\section{Twisted orthogonality}

We know from previous section that the cube has a certain quantum
symmetry algebra. This is by definition a quotient of $A_s(2^n)$,
the relations being those expressing edge preservation, in a
metric space or graph-theoretic sense.

Moreover, we have a preference for the cube regarded as a graph,
because this is a familiar object, namely the Cayley graph of the
group $L_n=Z_2^n$.

The quantum symmetry algebras of Cayley graphs are not computed in
general. However, for the group $L_n$, the study leads to the
following algebra:

\begin{definition}
$C(O_n^{-1})$ is the quotient of $A_o(n)$ by the following relations:
\begin{enumerate}
\item  $u_{ij}u_{ik}=-u_{ik}u_{ij}$, $u_{ji}u_{ki} =
-u_{ki}u_{ji}$,  for $i\neq j$. \item $u_{ij}u_{kl}=u_{kl}u_{ij}$
for $i\neq k$, $j\neq l$.
\end{enumerate}
\end{definition}

The relations in the statement are the skew-commutation relations
for the quantum group $GL_n^{-1}$. They define a Hopf ideal, so we
have indeed a Hopf algebra.

We should mention that this is not exactly the algebra of continuous fuctions on the quantum group $O_n^{-1}$, in the sense of classical texts like \cite{dr}, \cite{ma}, \cite{rft}. However, our definition is natural as well, and is further justified by the following result:

\begin{theorem}
$C(O_n^{-1})$ is the quantum symmetry algebra of $K_n$.
\end{theorem}

\begin{proof}
We have to show first that the quantum group $O_n^{-1}$ acts on
the cube. This means that the associated Hopf algebra
$C(O_n^{-1})$ coacts on the cube, in the sense explained in the
beginning of this section. Now recall from \cite{bb1} that this is
the same as asking for the existence of a coaction map:
$$\alpha :C(V_n)\to C(V_n)\otimes C(O_n^{-1})$$

Here $V_n$ is the set of vertices of $K_n$, and the coaction
$\alpha$ should be thought of as being the functional analytic
transpose of an action map $a:(v,g)\to g(v)$.

With notations from previous proof, we have to construct a map as
follows:
$$\alpha :C^*(L_n)\to C^*(L_n)\otimes C(O_n^{-1})$$

We claim that the following formula provides the answer:
$$\alpha(\tau_i)=\sum_j \tau_j
\otimes u_{ji}$$

Indeed, the elements $\tau_i$ generate the algebra on the left, with the relations $\tau_i^*=\tau_i$, $\tau_i^2=1$ and $\tau_i\tau_j=\tau_j\tau_i$. These relations are satisfied as well by the elements $\alpha(\tau_i)$, so we have a morphism of algebras. Now since $u$ is a corepresentation, $\alpha$ is coassociative and counital. It remains to check that $\alpha$ preserves the adjacency matrix of $K_n$.

It is routine to check that for $i_1 \not= i_2 \neq \ldots \neq
i_l$ we have:
$$\alpha(\tau_{i_1}\ldots\tau_{i_l})
=\sum_{j_1\neq\ldots\neq j_l}\tau_{j_1}\ldots\tau_{j_l} \otimes
u_{j_1i_1}\ldots u_{j_li_l}$$

In terms of eigenspaces $E_s$ of the adjacency matrix, this gives:
$$\alpha(E_s)\subset E_s\otimes C(O_n^{-1})$$

It follows that $\alpha$ preserves the adjacency matrix of $K_n$,
so it is a coaction on $K_n$ as claimed. See \cite{bb1} for
unexplained terminology in this proof.

So far, we know that $C(O_n^{-1})$ coacts on the cube. We have to show that the canonical map $B\to C(O_n^{-1})$ is an isomorphism, where $B=A_{K_n}$ is the quantum symmetry algebra of the cube. For this purpose, consider the universal coaction on the cube:
$$\beta:C^*(L_n)\to C^*(L_n)\otimes B$$

By eigenspace preservation we get elements $x_{ij}$ such that:
$$\beta(\tau_i)=\sum_j\tau_j\otimes x_{ji}$$

By applying $\beta$ to the relation $\tau_i\tau_j=\tau_j\tau_i$ we
get $x^tx=1$, so the matrix $x=(x_{ij})$ is orthogonal. By applying $\beta$
to the relation $\tau_i^2=1$ we get:
$$1\otimes\sum_kx_{ki}^2+\sum_{k<l}\tau_k\tau_l
\otimes(x_{ki}x_{li}+x_{li}x_{ki})=1\otimes 1$$

This gives $x_{ki}x_{li}=-x_{li}x_{ki}$ for $i\neq j$, $k\not=l$,
and by using the antipode we get $x_{ik}x_{il}=-x_{il}x_{ik}$ for
$k \not = l$. Also, by applying $\beta$ to
$\tau_i\tau_j=\tau_j\tau_i$ with $i\neq j$ we get:
$$\sum_{k<l} \tau_k\tau_l
\otimes(x_{ki}x_{lj}+x_{li}x_{kj})= \sum_{k<l}\tau_k \tau_l\otimes
(x_{kj}x_{li}+x_{lj}x_{ki})$$

It follows that for $i\not=j$ and $k\not=l$, we have:
$$x_{ki}x_{lj}+x_{li}x_{kj}=x_{kj}x_{li}+x_{lj}x_{ki}$$

In other words, we have $[x_{ki},x_{lj}] = [x_{kj},x_{li}]$. By
using the antipode we get $[x_{jl},x_{ik}]=[x_{il},x_{jk}]$. Now
by combining these relations we get:
$$[x_{il},x_{jk}]=[x_{ik},x_{jl}]=[x_{jk},x_{il}]=-[x_{il},x_{jk}]$$

This gives $[x_{il},x_{jk}]=0$, so the elements $x_{ij}$ satisfy
the relations for $C(O_n^{-1})$. This gives a map $C(O_n^{-1})\to B$, which is inverse to the canonical map $B\to C(O_n^{-1})$.
\end{proof}

We discuss now representations of $O_n^{-1}$. There are many techniques that can be applied, see e.g. \cite{ma}. We present here a simplified proof of the tensor
equivalence between the categories of finite dimensional representations of $O_n$ and $O_n^{-1}$, by using Morita
equivalence type techniques from \cite{bi1}, \cite{bdv},
\cite{sc}.

\begin{theorem}
The category of corepresentations of $C(O_n^{-1})$ is tensor equivalent
to the category of representations of $O_n$.
\end{theorem}

\begin{proof}
Let $A$ be the algebra of polynomial functions on $O_n$, the
standard matrix of generators being denoted $a_{ij}$. This is a dense
Hopf $*$-subalgebra of $C(O_n)$.

Let also $B$ be the universal $*$-algebra defined by the relations
of $C(O_n^{-1})$, with generators denoted $b_{ij}$. This embeds as
a dense Hopf $*$-algebra of $C(O_n^{-1})$.

The two categories in the statement are those of comodules over
$A$ and $B$. By \cite{sc}, what we have to do is to construct an
$A-B$ Hopf bi-Galois extension.

Let $Z$ be the universal algebra generated by the entries $z_{ij}$
of an orthogonal $n \times n$ matrix, subject to the following
relations:
$$z_{ij}z_{kl}=- z_{kl}z_{ij}\ {\rm for} \ j \neq l, \quad
z_{ij}z_{kj}=z_{kj}z_{ij}$$

We can make $Z$ into an $A-B$ bicomodule algebra, by constructing
coactions $\alpha:Z\to A\otimes Z$ and $\beta:Z\to Z\otimes B$ as
follows:
\begin{eqnarray*}
\alpha(z_{ij})&=&\sum_k a_{ik} \otimes z_{kj}\cr
\beta(z_{ij})&=&\sum_k z_{ik} \otimes b_{kj}
\end{eqnarray*}

We have to check first that $Z$ is nonzero. Consider the group
$L_n=<\tau_1,\ldots,\tau_n>$. We have a Hopf algebra morphism
$\pi:A\to C^*(L_n)$, given by $a_{ij}\to\delta_{ij}\tau_i$. We
define a bicharacter of $L_n$ as follows:
$$\sigma(\tau_i,\tau_j)=
\begin{cases}
-1\mbox{ if $i<j$}\\
+1\mbox{ if $i\geq j$}
\end{cases}$$

Since $\rho=\sigma(\pi\otimes\pi)$ is a Hopf 2-cocycle on $A$, we
can construct as in \cite{bi1} the right twisted algebra $A_\rho$.
The map $Z\to A_\rho$ given by $z_{ij}\to a_{ij}$ shows that $Z$
is nonzero.

Consider now the universal algebra $T$ generated by entries of an
orthogonal $n \times n$ matrix $t$, subject to the following
relations:
$$t_{ij}t_{kl}=-t_{kl}t_{ij}\ {\rm for} \ i \neq k, \quad
t_{ij}t_{il}=t_{il}t_{ij}$$

We define now algebra maps $\gamma:A\to Z\otimes T$ and
$\delta:B\to T\otimes Z$ as follows:
\begin{eqnarray*}
\gamma(a_{ij})&=&\sum_k z_{ik}\otimes t_{kj}\cr
\delta(b_{ij})&=&\sum_k t_{ik} \otimes z_{kj}
\end{eqnarray*}

Together with the map $s:T\to Z$ given by $t_{ij}\to z_{ji}$,
these make $(A,B,Z,T)$ into a Hopf-Galois system in the sense of
\cite{bi1}. Thus $Z$ is an $A-B$ bi-Galois object, and by
\cite{sc} this gives the desired equivalence of categories.
Moreover, our proof extends in a straightforward manner to the
$C^*$-algebraic setting of \cite{bdv}, so the tensor category
equivalence we have constructed is a tensor $C^*$-category
equivalence.
\end{proof}

The above result shows that integration over $O_n^{-1}$ appears as
a kind of twisting of integration over $O_n$. In particular, the
integral geometry of $O_n^{-1}$ is basically governed by
independence, so we don't have freeness.

\section{Sudoku unitaries}

We construct in this section the free version of $H_n$. The idea
is to use the fact that $H_n$ is the symmetry group of $I_n$, the
graph formed by $n$ segments.

\begin{definition}
A sudoku unitary is a magic unitary of the form
$$w=\begin{pmatrix}a&b\cr b&a\end{pmatrix}$$
where both $a,b$ are square matrices.
\end{definition}

In this definition the matrices $a,b$ are of course of same size. The terminology comes from a vague analogy with sudoku squares. 

As a first example, any $2\times 2$ magic unitary is a sudoku unitary:
$$w=\begin{pmatrix}p&1-p\cr 1-p&p\end{pmatrix}$$

Here is an example of $4\times 4$ sudoku unitary, with $p,q$ being
projections:
$$w=\begin{pmatrix}p&0&1-p&0\cr 0&q&0&1-q\cr 1-p&0&p&0\cr
0&1-q&0&q\end{pmatrix}$$

The third example comes from the hyperoctahedral group. When
regarding $H_n$ as symmetry group of $I_n$, hence as subgroup of
$S_{2n}\subset O_{2n}$, the matrix coordinates form a sudoku
unitary. Indeed, the characteristic functions $a_{ij},b_{ij}$ are
those of group elements mapping $I^i\to I^j$ by fixing
orientations, respectively by changing them.

\begin{definition}
$A_h(n)$ is the universal algebra generated by the entries of a
$2n\times 2n$ sudoku unitary.
\end{definition}

It follows from the above discussion that we have $A_h(1)=C(Z_2)$,
and that the algebra $A_h(2)$ is not commutative, and infinite
dimensional. The noncommutativity and infinite dimensionality of
$A_h(n)$ hold in fact for any $n\geq 2$, because we have surjective morphisms $A_h(n+1)\to A_h(n)$. These are given by the fact that if $w=(^a_b{\ }^b_a)$ is sudoku, then the following matrix is sudoku as well:
$$w'=\begin{pmatrix}a&&b&\cr &1&&0\cr b&&a&\cr &0&&1\end{pmatrix}$$

Here the blank spaces are filled with row or column vectors consisting of zeroes.

A first result motivating our definition is as follows.

\begin{proposition}
$A_h(n)$ is the quantum symmetry algebra of $I_n$.
\end{proposition}

\begin{proof}
The quantum symmetry algebra of $I_n$ is by definition the
quotient of $A_s(2n)$ by the relations coming from commutation of
the magic unitary matrix, say $w$, with the adjacency matrix of
$I_n$. This adjacency matrix is given by:
$$e=\begin{pmatrix}0&1_n\cr 1_n&0\end{pmatrix}$$

We can write as well $w$ as a $2\times 2$ matrix, with $n\times n$
matrix entries:
$$w=\begin{pmatrix}a&b\cr c&d\end{pmatrix}$$

The relation $ew=we$ is equivalent to $a=d,b=c$, and this gives
the result.
\end{proof}

We proceed now with the study of $A_h(n)$, in analogy with what we
know about $H_n$. Our first goal is to extend the wreath product decomposition $H_n=Z_2\wr S_n$. We use the notion of free wreath product, introduced in \cite{bi2}. 

Consider two quantum permutation
algebras $(A,u)$ and $(B,v)$. This means that we have surjective Hopf algebra morphisms $A_s(n)\to A$ and $A_s(m)\to B$
for certain numbers $n,m$, and the fundamental corepresentations of $A,B$ are denoted $u,v$.

\begin{definition}
The free wreath product of $(A,u)$ and $(B,v)$ is given by
$$A*_wB=(A^{*n}*B)/<[u_{ij}^{(a)},v_{ab}]=0>$$
where $n$ is the size of $v$, and has magic unitary matrix
$w_{ia,jb}=u_{ij}^{(a)}v_{ab}$.
\end{definition}

This definition is justified by formulae of the following type,
where $G,A$ denote classical symmetry groups, respectively quantum
symmetry algebras:
\begin{eqnarray*}
G(X*Y)&=&G(X)\ \wr\  G(Y)\cr A(X*Y)&=&A(X)*_wA(Y)\end{eqnarray*}

There are several such formulae, depending on the types of graphs
and free products considered. See \cite{bb1},
\cite{bi2}. The formula we are interested in is:
\begin{eqnarray*}
G(I\ldots I)&=&G(I)\ \wr \ G(\circ\ldots\circ)\cr A(I\ldots
I)&=&A(I)*_wA(\circ\ldots\circ)\end{eqnarray*}

Here $I$ is a segment, $\circ$ is a point, and the dots mean
$n$-fold disjoint union. We recognize at left the graph $I_n$. The
upper formula is nothing but the isomorphism $H_n=Z_2\wr
S_n$ we started with. As for the lower formula, this is what we
need:

\begin{theorem}
We have $A_h(n)=C(Z_2)*_wA_s(n)$.
\end{theorem}

This follows from the above discussion, see \cite{bi2}, \cite{bb1}. As a first application, we can make a comparison between $A_h(n)$
and $C(O_n^{-1})$. Recall that both algebras appear as
noncommutative versions of $C(H_n)$, in the sense that their
maximal commutative quotient if $C(H_n)$. The problem is whether these algebras are isomorphic or not.

\begin{proposition}
We have the following isomorphisms:
\begin{enumerate}
\item $A_h(2)=C(Z_2)*_wC(Z_2)$.
\item $A_h(3)=C(Z_2)*_wC(S_3)$.
\item $C(O_2^{-1})=C(Z_2)*_wC(Z_2)$.
\item $C(O_3^{-1})=C(Z_2)\otimes A_s(4)$.
\end{enumerate}
\end{proposition}

\begin{proof}
The first two assertions follow from the above theorem and from the exceptional isomorphisms $A_s(2)=C(Z_2)$ and $A_s(3)=C(S_3)$.

The last two assertions follow by using the exceptional isomorphisms $K_2=I^c$ and $K_3=I\times T$, where $I$ is the segment and $T$ is the tetrahedron. See \cite{bb2} for details.
\end{proof}

This result shows that $A_h(n)$ and $C(O_n^{-1})$
are equal for $n=2$, and different for $n=3$. Actually these
algebras are different for any $n\geq 3$, as one can see by
comparing representation theory invariants: the fusion rules for
$C(O_n^{-1})$ are the same as those for $O_n$, and in particular
they are commutative, while for $A_h(n)$ with $n\geq 3$ we will
see that the fusion rules are not commutative.

\section{Cubic unitaries}

In this section we present an alternative approach to $A_h(n)$.
The idea is to use the fact that $H_n$ is a subgroup of $O_n$,
hence the matrix of coordinates has size $n$. The relevant
condition on the entries of this matrix is as follows.

\begin{definition}
A cubic unitary is an orthogonal matrix $u$ satisfying
$$u_{ij}u_{ik}=u_{ji}u_{ki}=0$$
for $j\neq k$. In other words, on rows and columns, distinct
entries satisfy $ab=ba=0$.
\end{definition}

The basic example comes indeed from the hyperoctahedral group
$H_n\subset O_n$. In terms of the action on the space
$I_n=\{I^1,\ldots ,I^n\}$ formed by the $[-1,1]$ segments on each
coordinate axis, the matrix coordinates are given by:
$$u_{ij}(g)=\begin{cases}\ 0\ \mbox{ if $g(I^j)\neq I^i$}\\
+1\mbox{ if $g(I^j)=I^i$}\\
-1\mbox{ if $g(I^j)=I^i$ returned}
\end{cases}$$

Now when $i$ is fixed and $j$ varies, or vice versa, the supports
of $u_{ij}$ form partitions of $H_n$. Thus the matrix of
coordinates $u$ is cubic in the above sense.

\begin{theorem}
$A_h(n)$ is the universal algebra generated by the entries of a
$n\times n$ cubic unitary matrix.
\end{theorem}

\begin{proof}
Consider the universal algebra in the statement:
$$A_c(n)=C^*\left(u_{ij}\ \Big\vert\  u=n\times n\mbox{ cubic
unitary}\right)$$

We have to compare it with the following universal algebra:
\begin{eqnarray*}
A_h(n)&=&C^*\left(w_{ij}\ \Big\vert\  w=2n\times 2n\mbox{ sudoku
unitary}\right)\cr &=&C^*\left(a_{ij},b_{ij}\
\Big\vert\begin{pmatrix}a&b\cr b&a\end{pmatrix}=2n\times 2n\mbox{
magic unitary}\right)
\end{eqnarray*}

We construct first the arrow $A_c(n)\to A_h(n)$. Consider the
matrix $a-b$, having as entries the elements $a_{ij}-b_{ij}$. The
elements $a_{ij},b_{ij}$ being self-adoint, their difference is
self-adjoint as well. Thus $a-b$ is a matrix of self-adjoint
elements. We have the following formula for products on columns of
$a-b$:
\begin{eqnarray*}
(a-b)_{ik}(a-b)_{jk}
&=&a_{ik}a_{jk}-a_{ik}b_{jk}-b_{ik}a_{jk}+b_{ik}b_{jk}\cr &=&
\begin{cases}
0-0-0+0\mbox{ for }i\neq j\\
a_{ik}-0-0+b_{ik}\mbox{ for }i=j \end{cases}\cr &=&\begin{cases}
0\mbox{ for }i\neq j\\
a_{ik}+b_{ik}\mbox{ for }i=j
\end{cases}
\end{eqnarray*}

In the $i=j$ case it follows from the sudoku condition that the
elements $a_{ik}+b_{ik}$ sum up to $1$. Thus the columns of $a-b$
are orthogonal. A similar computation works for rows, and we
conclude that $a-b$ is an orthogonal matrix.

Now by using the $i\neq j$ computation, along with its row
analogue, we conclude that $a-b$ is cubic. Thus we can define our
morphism by the following formula:
$$\varphi(u_{ij})=a_{ij}-b_{ij}$$

We construct now the inverse morphism. Consider the following
elements:
\begin{eqnarray*}
\alpha_{ij}&=&\frac{u_{ij}^2+u_{ij}}{2}\cr
\beta_{ij}&=&\frac{u_{ij}^2-u_{ij}}{2}
\end{eqnarray*}

It follows from the cubic condition that these elements are
projections, and that the following matrix is a sudoku unitary:
$$M=\begin{pmatrix}(\alpha_{ij})&(\beta_{ij})\cr
(\beta_{ij})&(\alpha_{ij})\end{pmatrix}$$

Thus we can define our morphism by the following formula:
\begin{eqnarray*}
\psi(a_{ij})&=&\frac{u_{ij}^2+u_{ij}}{2}\cr
\psi(b_{ij})&=&\frac{u_{ij}^2-u_{ij}}{2}
\end{eqnarray*}

We check now the fact that $\psi,\varphi$ are indeed inverse
morphisms:
\begin{eqnarray*}
\psi\varphi(u_{ij})&=&\psi(a_{ij}-b_{ij})\cr
&=&\frac{u_{ij}^2+u_{ij}}{2}-\frac{u_{ij}^2-u_{ij}}{2}\cr
&=&u_{ij}
\end{eqnarray*}

As for the other composition, we have the following computation:
\begin{eqnarray*}
\varphi\psi(a_{ij})&=&\varphi\left(\frac{u_{ij}^2+u_{ij}}{2}\right)\cr
&=&\frac{(a_{ij}-b_{ij})^2+(a_{ij}-b_{ij})}{2}\cr
&=&a_{ij}
\end{eqnarray*}

A similar computation gives $\varphi\psi(b_{ij})=b_{ij}$, which
completes the proof.
\end{proof}

We are interested now in quantum analogues of the embeddings $S_n\subset
H_n\subset O_n$. It is convenient to state the
corresponding result in the following way.

\begin{proposition}
We have the following commutative diagram:
$$\begin{matrix}
A_o(n)&\to&A_h(n)&\to&A_s(n)\cr\cr
\downarrow&&\downarrow&&\downarrow\cr\cr
C(O_n)&\to&C(H_n)&\to&C(S_n)
\end{matrix}$$

Moreover, both horizontal compositions are the canonical quotient
maps.
\end{proposition}

\begin{proof}
The upper left arrow $A_o(n)\to A_h(n)$ is by definition the one
coming from previous theorem. We construct now the arrow
$A_h(n)\to A_s(n)$. With sudoku notations for $A_h(n)$, consider
the ideal $I$ generated by the relations $b_{ij}=0$. We have:
\begin{eqnarray*}
A_h(n)/I&=&C^*\left(a_{ij}\ \Big\vert\begin{pmatrix}a&0\cr
0&a\end{pmatrix}=2n\times 2n\mbox{ magic
unitary}\right)\cr&=&C^*\left(a_{ij}\ \Big\vert\ a=n\times n\mbox{
magic unitary}\right)\end{eqnarray*}

The algebra on the right is $A_s(n)$, and this gives a map
$A_h(n)\to A_s(n)$. The fact that both squares of the diagram
commute is clear from definitions.

It remains to check the last assertion. With the above notations, the upper
horizontal composition is given by $u_{ij}\to a_{ij}-b_{ij}\to
a_{ij}$, so it is the canonical quotient map $A_o(n)\to A_s(n)$.
As for the lower horizontal composition, the fact that this is the
transpose of the canonical embedding $S_n\subset O_n$ follows from definitions.
\end{proof}

One can prove that the quotient map $A_h(n)\to A_s(n)$ has a
canonical section, which is a quantum analogue for the transpose
of the canonical quotient map $H_n\to S_n$.

The general idea is that the quantum analogues of various
properties of $H_n$ are shared by the quantum groups $O_n^{-1}$
and $H_n^+$, the latter being the spectrum of $A_h(n)$.

\section{Representation theory}

For the rest of this paper, we make the assumption $n\geq 4$.

In this section we discuss representation theory problems. The
idea is that the presentation relations for various quantum
algebras have graphical interpretations, which lead via Tannakian
duality to Temperley-Lieb type diagrams.

The fact that free wreath products of type $A_h(n)$ correspond to
the Fuss-Catalan algebra of Bisch and Jones \cite{bj} is known from some time. See \cite{bb1} for a general
statement in this sense, and for the complete story and
state-of-art of the subject.

In this paper we have some special requirements in this sense:
\begin{enumerate}
\item In order to get the good integration formulae, the
fundamental corepresentation here must be the $n\times n$ cubic
unitary. This is a bit different from what can be found in the
litterature, where the $2n\times 2n$ sudoku unitary is used. \item
In order to use various free probability formulae, we need results
in terms of colored non-crossing partitions, rather than in terms
of colored Temperley-Lieb diagrams (or Fuss-Catalan diagrams) used
in the litterature.\item The whole thing is to be extended in next
sections, for Hopf algebras satisfying $A_o(n)\to A\to A_s(n)$. In
other words, the results for $A_o(n),A_s(n)$ must be explained as
well in detail, at some point of the paper.
\end{enumerate}

Summarizing, we need a uniform treatment of the algebras $A_x(n)$
with $x=o,h,s$, with the fundamental corepresentation being the
$n\times n$ one, and in terms of both diagrams and partitions. We
use the group subscripts $o,h,s$ for various types of diagrams.

\begin{definition}
The set of Temperley-Lieb diagrams is
$$TL_s=\left\{\begin{matrix}\cdot\,\cdot\,\cdot\,\,\cdot& \leftarrow &
    2l\,\,\mbox{points}\cr W & \leftarrow &
    l+k\mbox{ strings}\cr \cdot\,\cdot\,\cdot\,\cdot\,\cdot\,\,\cdot& \leftarrow &
    2k\,\,\mbox{points}\end{matrix}\right\}$$
where strings join pairs of points, and do not cross.
\end{definition}

Given a Temperley-Lieb diagram, we can color both rows of points $abbaabba..$ with the pattern
$abba$ repeated as many time as needed, and with an extra $ab$
block at the end, if needed. If strings can be colored as well, we say that the diagram is colored.

\begin{definition}
The set of colored Temperley-Lieb diagrams is
$$TL_h=\left\{\begin{matrix}\cdot\,\cdot\,\cdot\,\,\cdot& \leftarrow &
    2l\,\,\mbox{points}\cr W & \leftarrow &
    l+k\mbox{ colored strings}\cr \cdot\,\cdot\,\cdot\,\cdot\,\cdot\,\,\cdot& \leftarrow &
    2k\,\,\mbox{points}\end{matrix}\right\}$$
where both rows of points are colored $abbaabba..$.
\end{definition}

Given a Temperley-Lieb diagram, we can collapse pairs of consecutive neighbors, in both rows. If strings can be collapsed as well, we say that the diagram is double. Observe that a double Temperley-Lieb diagram is colored.

\begin{definition}
The set of double Temperley-Lieb diagrams is
$$TL_o=\left\{\begin{matrix}\cdot\,\cdot\,\cdot\,\,\cdot& \leftarrow &
    2l\,\,\mbox{points}\cr W & \leftarrow &
    (l+k)/2\mbox{ double strings}\cr \cdot\,\cdot\,\cdot\,\cdot\,\cdot\,\,\cdot& \leftarrow &
    2k\,\,\mbox{points}\end{matrix}\right\}$$
where points and strings are taken, as usual, up to isotopy.
\end{definition}

Probably most illustrating here is the following table, for small number of points. The Temperley-Lieb diagrams between $2l$ points and $2k$ points are called $l\to k$ diagrams.
\begin{center}
\begin{tabular}[t]{|l|l|l|l|l|}
\hline &$TL_o$&$TL_h$&$TL_s$&$\# TL_s$\\
\hline $0\to 1$&&&$\cap$&1\\
\hline $0\to 2$&$\Cap$&$\Cap$&$\cap\cap$ , $\Cap$&2\\
\hline $1\to 1$&$||$&$||$&$||$ , ${\ }^\cup_\cap$&2\\
\hline $1\to 2$&&&$||{\,}_\cap$ , ${\,}_\cap ||$ , $|{\,}_\cap|$ , ${\ }^{\;\cup}_{\cap\cap}$ , ${\,}^\cup_\Cap$&5\\
\hline $2\to 2$&$||||$ , $^\Cup_\Cap$&$||||$ , $^\Cup_\Cap$ ,
$|\!{\;}^\cup_\cap|$&$||||$ , $^\Cup_\Cap$ , $|\!{\;}^\cup_\cap|$
, $^{\;\Cup}_{\cap\cap}$ ,
${\,}^\cup\!/\!\!{\;}_\cap$ , $\dots$&14\\
\hline
\end{tabular}
\end{center}
\bigskip

We recognize at right the Catalan numbers. The number of elements
in the left column is also a Catalan number. As for the middle
column, the diagrams here are counted by Fuss-Catalan numbers.
This will be discussed in detail later on.

To any number $n$ and any subscript $x\in\{o,h,s\}$ we associate a
tensor category $TL_x(n)$, in the following way:
\begin{enumerate}
\item The objects are the positive integers: $0,1,2,\ldots$ \item
The arrows $l\to k$ are linear combinations of $l\to k$ diagrams in $TL_x$. \item The composition is by vertical concatenation,
with the rule $\bigcirc =n$. \item The tensor product is by
horizontal concatenation. \item The involution is by upside-down
turning.
\end{enumerate}

Temperley-Lieb diagrams act on tensors according to the following
formula, where the middle symbol is $1$ if all strings of $p$ join
pairs of equal indices, and is $0$ if not:
$$p(e_{i_1}\otimes\ldots\otimes e_{i_l})
=\sum_{j_1\ldots j_k}\begin{pmatrix} i_1\,i_1\ldots i_l\,i_l\cr
p\cr j_1\,j_1\ldots
j_k\,j_k\end{pmatrix}e_{j_1}\otimes\ldots\otimes e_{j_k}$$

By using the trace one can see that different diagrams produce
different linear maps, and this action makes $TL_x(n)$ a subcategory
of the category of Hilbert spaces.

\begin{theorem}
We have the Tannakian duality diagram
$$\begin{matrix}
A_o(n)&\to&A_h(n)&\to&A_s(n)\cr\cr
\updownarrow&&\updownarrow&&\updownarrow\cr\cr
TL_o(n)&\subset&TL_h(n)&\subset&TL_s(n)
\end{matrix}$$
where all horizontal arrows are the canonical ones.
\end{theorem}

The proof uses Woronowicz's Tannakian duality in \cite{wo2}. The
idea is that the presentations of the above algebras translate
into presentations of associated tensor categories, which in turn
lead to diagrams. Consider a unitary matrix $u$.
\begin{enumerate}
\item  The fact that $u$ is orthogonal is equivalent to $V\in
Hom(1,u^{\otimes 2})$, where $V$ is the canonical vector in
$\mathbb C^n\otimes\mathbb C^n$. But the relations satisfied by
$V$ in a categorical sense are those satisfied by $v=\Cap$, which
generates $TL_o$. See \cite{bc1}.\item The fact that $u$ is cubic
is equivalent to $V\in Hom(1,u^{\otimes 2})$ and $E\in
End(u^{\otimes 2})$, where $E$ is the canonical projection. But
the relations satisfied by $V,E$ are those satisfied by $v=\Cap$
and $e=|\!{\;}^\cup_\cap|$, which generate $TL_h$. See
\cite{ba2}.\item The fact that $u$ is magic is equivalent to $M\in
Hom(u^{\otimes 2},u)$ and $U\in Hom(1,u)$, where $M,U$ are the
multiplication and unit of $\mathbb C^n$. But $M,U$ satisfy the same relations as $m=|\cup|$ and $u=\cap$,
which generate $TL_s$. See \cite{bc2}.
\end{enumerate}

These three remarks give the three vertical arrows, standing for Tannakian dualities. As for horizontal compatibility, this is clear from definitions.

\section{The Weingarten formula}

We can proceed now with the Weingarten general algorithm. Our first task is to translate the Tannakian result into a statement regarding fixed points.

\begin{definition}
An arch is a Temperley-Lieb diagram having no upper points.
\end{definition}

The set of double, colored and usual arches are denoted $AR_x$
with $x=o,h,s$. We have the following table, for small values of
$k$.
\begin{center}
\begin{tabular}[t]{|l|l|l|l|l|}
\hline &$AR_o$&$AR_h$&$AR_s$&$\# AR_s$\\
\hline $1$&&&$\cap$&1\\
\hline $2$&$\Cap$&$\Cap$&$\Cap$ , $\cap\cap$&2\\
\hline $3$&&&$\cap\Cap$ , $\Cap\cap$ , $\cap\cap\cap$ , $|^{{\
}^{\!\!\!\!-\!\!-\!\!-\!\!-\!\!}}| \hskip-5.6mm{\ }_{\Cap}$ \ \ ,
$|^{{\ }^{\!\!\!\!-\!\!-\!\!-\!\!-\!\!-\!\!}}|
\hskip-7mm{\ }_{{\cap\cap}}$\ \ &5\\
\hline $4$&$\Cap\Cap$ , $|^{{\ }^{\!\!\!\!-\!\!-\!\!-\!\!-\!\!}}|
\hskip-6.6mm{\ }_{{\ }_\Cap}\hskip-3.8mm{\ }_{\bigcap}$ \ \
&$\Cap\Cap$ , $|^{{\ }^{\!\!\!\!-\!\!-\!\!-\!\!-\!\!}}|
\hskip-6.6mm{\ }_{{\ }_\Cap}\hskip-3.8mm{\ }_{\bigcap}$ \ \ ,
$|^{{\ }^{\!\!\!\!-\!\!-\!\!-\!\!-\!\!-\!\!-\!\!}}| \hskip-8.6mm{\
}_{{\cap\cap\cap}}$ \  &$\Cap\Cap$ , $|^{{\
}^{\!\!\!\!-\!\!-\!\!-\!\!-\!\!}}| \hskip-6.6mm{\ }_{{\
}_\Cap}\hskip-3.8mm{\ }_{\bigcap}$ \ \ , $|^{{\
}^{\!\!\!\!-\!\!-\!\!-\!\!-\!\!-\!\!-\!\!}}| \hskip-8.6mm{\
}_{{\cap\cap\cap}}$\ \  , $\cap\Cap\cap$ , $\dots$&14\\
\hline
\end{tabular}
\end{center}
\bigskip

The Tannakian duality result in previous section has the following
reformulation.

\begin{proposition}
The space of fixed points of $u^{\otimes k}$ is as follows.
\begin{enumerate}
\item For $A_o(n)$ this is the space spanned by $AR_o(k)$. \item
For $A_h(n)$ this is the space spanned by $AR_h(k)$. \item For
$A_s(n)$ this is the space spanned by $AR_s(k)$.
\end{enumerate}
\end{proposition}

The arches can be regarded as non-crossing partitions, by
collapsing pairs of consecutive neighbors. As an example, we have
the following table:
\begin{center}
\begin{tabular}[t]{|l|l|l|l|l|l|l|l|l|}
\hline arch&$\cap$&$\Cap$&$\cap\cap$& $\Cap\cap$&$\cap\Cap$&
$|^{{\ }^{\!\!\!\!-\!\!-\!\!-\!\!-\!\!-\!\!-\!\!}}_
{\ {\cap\cap}}|$&$\cap\cap\cap$\\
\hline partition&$|$&$\cap$&$|\ |$& $\cap\  |$&$|\ \cap$&$\ \
\sqcap\hskip-0.6mm\sqcap$&
$|\ \ |\ \ |$\\
 \hline
\end{tabular}
\end{center}
\bigskip

This construction probably deserves more explanations. First, we
regard non-crossing partitions as diagrams, in the obvious way.
Here is for instance the diagram corresponding to the partition
$\{1,2,3,4,5,6\}=\{1,2,5\}\cup\{3,4\}\cup\{6\}$:
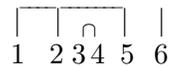
$$\begin{matrix}|^{{\ }^{\!\!\!\!-\!\!-\!\!-\!\!-\!\!}}_
{\ }|^{{\ }^{\!\!\!\!-\!\!-\!\!-\!\!-\!\!-\!\!-\!\!-\!\!}}_{\ \
\cap}|\ \ \ |\cr
 1\ \ \, 2\; 3\,4\ \,5\ \ 6\end{matrix}$$

 And here is an example of fattening operation, producing
 arches from partitions:
 $$|\ \ \sqcap\hskip-1.4mm\sqcap\ \ \cap
 \ \longrightarrow\
 ||\ \ ||^{=}||^{=}||\ \ ||^{=}||
 \ \longrightarrow\
 \cap\ |^{{\ }^{\!\!\!\!-\!\!-\!\!-\!\!-\!\!-\!\!-\!\!}}_{\ \cap\cap}|\ \Cap$$

The sets of partitions coming from double arches, colored arches and usual arches are denoted $NC_x(k)$, with $x=o,h,s$. Here $k$ is the number of points of the set which is partitioned, equal to half of the number of lower points of the arch.

The above result has the following reformulation.

\begin{proposition}
The space of fixed points of $u^{\otimes k}$ is as follows:
\begin{enumerate}
\item For $A_o(n)$ this is the space spanned by $NC_o(k)$.
\item For $A_h(n)$ this is the space spanned by $NC_h(k)$.
\item For $A_s(n)$ this is the space spanned by $NC_s(k)$.
\end{enumerate}
\end{proposition}

It is useful to give an intrinsic description of all these
notions. Going from arches to non-crossing partitions leads to the
following sets: \begin{enumerate} \item $NC_o(k)$ is the set of
non-crossing partitions into pairs. \item $NC_h(k)$ is the set of non-crossing partitions into even parts. \item $NC_s(k)$ is the set of
non-crossing partitions.
\end{enumerate}

The elements of $NC_o(k)$ are also called non-crossing pair-partitions, and those of $NC_h(k)$, non-crossing colored partitions. This is because we have the following characterization of $NC_h(k)$, which makes the link with colored diagrams: when attaching alternative labels $ab,ba,ab,ba,\ldots$ to the legs of
the partition, the condition is that when travelling on sides of
strings, different labels do not meet. Here is an example:
$$\begin{matrix}|^{{\ }^{\!\!\!\!-\!\!-\!\!-\!\!-\!\!-\!\!-\!\!}}
|^{{\ }^{\!\!\!\!-\!\!-\!\!-\!\!-\!\!-\!\!-\!\!}}|^{{\
}^{\!\!\!\!-\!\!-\!\!-\!\!-\!\!-\!\!-\!\!}}| \cr a\,b\ \ \ b\,a\ \
\ a\,b\ \ \ b\,a\end{matrix}\ \ \
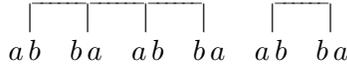
\begin{matrix}|^{{\ }^{\!\!\!\!-\!\!-\!\!-\!\!-\!\!-\!\!-\!\!}}|
\cr a\,b\ \ \ b\,a\end{matrix}$$

Probably most illustrating here is the following table, for small
values of $k$.
\begin{center}
\begin{tabular}[t]{|l|l|l|l|l|}
\hline &$NC_o$&$NC_h$&$NC_s$&$\# NC_s$\\
\hline $1$&&&$|$&1\\
\hline $2$&$\cap$&$\cap$&$\cap$ , $|\ |$&2\\
\hline $3$&&&$|\cap$ , $\cap |$ , $|\ |\ |$ ,
$\bigcap\hskip-5.2mm{\ }_{{\ }^{|}}$ \ , $\sqcap\hskip-.6mm\sqcap$&5\\
\hline $4$&$\cap\cap$ , $\Cap$&$\cap\cap$ , $\Cap$ ,
$\sqcap\hskip-1.4mm\sqcap\hskip-1.4mm\sqcap$ &$\cap\cap$ , $\Cap$ ,
$\sqcap\hskip-1.4mm\sqcap\hskip-1.4mm\sqcap$ , $|\cap|$ , $\dots$&14\\
\hline
\end{tabular}
\end{center}
\bigskip

Now regarding fixed vectors, we can plug multi-indices
$i=(i_1\ldots i_{k})$ into partitions of $\{1,\ldots ,k\}$, then
define coupling numbers $\delta_{pi}$ and vectors $v_{p}$, in the
following way.

\begin{proposition}
The vector created by a partition is given by
$$v_p=\sum_{i}\delta_{pi}\,e_{i_1}\otimes\ldots\otimes e_{i_k}$$
where $\delta_{pi}=0$ if some block of $p$ contains two different indices of $i$, and $\delta_{pi}=0$ if not.
\end{proposition}

This is indeed a reformulation of the action of Temperley-Lieb
diagrams on tensors: the general formula makes arches create
tensors, and by using the correspondence between arches and
non-crossing partitions we get the above formula.

For a partition $p$ we denote by $|p|$ the number of its blocks.

For two partitions $p,q$, we define $p\vee q$ to be the partition
obtained as set-theoretic least upper bound of $p,q$. That is,
$p\vee q$ is such that two elements which are in the same block of
$p$, or in the same block of $q$, have to be in the same block of
$p\vee q$.

\begin{theorem}
We have the Weingarten formula
$$\int u_{i_{1}j_{1}}\ldots u_{i_{k}j_{k}}=\sum_{pq}
\delta_{pi}\delta_{qj}W_{kn}(p,q)$$ where $u$ is the fundamental
corepresentation of $A_x(n)$, and:
\begin{enumerate}
\item The diagrams are the partitions in $NC_x(k)$. \item The
couplings with indices are by plugging indices into partitions.
\item The Gram matrix is given by $G_{kn}(p,q)=n^{|p\vee
q|}$.\item The Weingarten matrix is given by
$W_{kn}=G_{kn}^{-1}$.\end{enumerate}
\end{theorem}

\begin{proof}
The fact that $G_{kn}$ is indeed the Gram matrix is checked as follows:
\begin{eqnarray*}
<v_p,v_q> &=&\left<\sum_{i}\delta_{pi}\,e_{i_1}\otimes
\ldots\otimes e_{i_k},\sum_{i}\delta_{qi}\,e_{i_1}\otimes
\ldots\otimes e_{i_k}\right>\cr&=&
\sum_i\delta_{pi}\delta_{qi}=\sum_i\delta_{p\vee q,i}=n^{|p\vee q|}
\end{eqnarray*}

The rest of the statement is clear from definitions and from
general theory in \cite{wo1}, as explained in section 1 for
classical groups. See \cite{bc1} or \cite{bc2} for details.
\end{proof}

\begin{theorem}
We have the moment formula
$$\int (u_{11}+\ldots +u_{ss})^{k}=Tr(G_{kn}^{-1}G_{ks})$$
where $u$ is the fundamental corepresentation of $A_x(n)$.
\end{theorem}

\begin{proof}
We have the following computation:
\begin{eqnarray*}
\int (u_{11}+\ldots +u_{ss})^{k}
&=&\sum_{i_1=1}^{s}\ldots\sum_{i_{k}=1}^s\sum_{pq} \delta_{pi}
\delta_{qi}G_{kn}^{-1}(p,q)\cr
&=&\sum_{pq}G_{kn}^{-1}(p,q)\sum_{i_1=1}^{s}\ldots\sum_{i_{k}=1}^s
\delta_{pi}\delta_{qi}\cr &=&\sum_{pq}G^{-1}_{kn}(p,q)G_{ks}(p,q)
\end{eqnarray*}

This gives the formula in the statement.
\end{proof}

\section{Numeric results}

We can compute now the low order moments of
$u_{11}+\ldots+u_{ss}$. For $A_o(n)$ and $A_s(n)$ this is done in
\cite{bc1}, \cite{bc2}. For $A_h(n)$ the situation is as follows. 

\begin{theorem}
For the algebra $A_h(n)$ we have the moment formulae
\begin{eqnarray*}
\int (u_{11}+\ldots +u_{ss})^2\ &=&\frac{s}{n}\cr \int
(u_{11}+\ldots +u_{ss})^4&=&\frac{s}{n}\cdot \frac{n+2s-3}
{n-1}\cr \int (u_{11}+\ldots +u_{ss})^6&=&\frac{s}{n}\cdot
\frac{P_s(n)}{(n-1)^2(n-2)}
\end{eqnarray*}
where $P_s(n)=n^3+(6s-10)n^2+(5s^2-30s+30)n-(8s^2-30s+24)$.
\end{theorem}

\begin{proof}
At $l=1$ there is just one diagram, namely $\cap$, having one
block. Thus the Gram matrix is $G_{2n}=(n)$, and the computation
of the second moment goes as follows:
$$Tr(G_{2n}^{-1}G_{2s})=Tr\left(
(n^{-1})(s)\right)=Tr(n^{-1}s)=n^{-1}s$$

At $l=2$ we have three diagrams. The $p\vee q$ diagrams are
computed as follows:
\begin{center}
\begin{tabular}[t]{|l|l|l|l|l|}
\hline &$\cap\cap$&$\Cap$&
$\sqcap\hskip-1.4mm\sqcap\hskip-1.4mm\sqcap$\\
\hline $\cap\cap$&$\cap\cap$&$\sqcap\hskip-1.4mm\sqcap\hskip-1.4mm\sqcap$&$\sqcap\hskip-1.4mm\sqcap\hskip-1.4mm\sqcap$\\
\hline $\Cap$&$\sqcap\hskip-1.4mm\sqcap\hskip-1.4mm\sqcap$&$\Cap$&$\sqcap\hskip-1.4mm\sqcap\hskip-1.4mm\sqcap$\\
\hline $\sqcap\hskip-1.4mm\sqcap\hskip-1.4mm\sqcap$&$\sqcap\hskip-1.4mm\sqcap\hskip-1.4mm\sqcap$&$\sqcap\hskip-1.4mm\sqcap\hskip-1.4mm\sqcap$&$\sqcap\hskip-1.4mm\sqcap\hskip-1.4mm\sqcap$\\
\hline
\end{tabular}
\end{center}
\bigskip

Thus the Gram and Weingarten matrices are:
$$G_{4n}=n\begin{pmatrix}n&1&1\cr 1&n&1\cr 1&1&1\end{pmatrix}
\hskip 1cm G^{-1}_{4n}=\frac{1}{n(n-1)}
\begin{pmatrix}1&0&-1\cr 0&1&-1\cr -1&-1&n+1\end{pmatrix}$$

This gives the formula of the fourth moment:
\begin{eqnarray*}Tr(G_{4n}^{-1}G_{4s})&=&\frac{s}{n(n-1)}
\ Tr\left( \begin{pmatrix}1&0&-1\cr 0&1&-1\cr
-1&-1&n+1\end{pmatrix}
\begin{pmatrix}s&1&1\cr 1&s&1\cr 1&1&1\end{pmatrix}\right)\cr
&=&\frac{s}{n(n-1)}\ Tr
\begin{pmatrix}s-1&*&*\cr *&s-1&*\cr *&*&n-1\end{pmatrix}\cr
&=&\frac{s(n+2s-3)}{n(n-1)}
\end{eqnarray*}

At $l=3$ we have $12$ diagrams, and after computing the
corresponding $p\vee q $ diagrams we get the following Gram
matrix, whose inverse will not be given here:
$$G_{6n}=n\left(\begin{array}{ccccrccccccccc}
n^2&n&n&n&n&1&n&1&n&1&1&1\\
n&n^2&1&n&1&n&1&1&1&n&1&1\\
n&1&n^2&1&n&n&1&1&1&1&n&1\\
n&n&1&n&1&1&1&1&1&1&1&1\\
n&1&n&1&n&1&1&1&1&1&1&1\\
1&n&n&1&1&n^2&n&n&1&n&n&1\\
n&1&1&1&1&n&n^2&n&n&1&1&1\\
1&1&1&1&1&n&n&n&1&1&1&1\\
n&1&1&1&1&1&n&1&n&1&1&1\\
1&n&1&1&1&n&1&1&1&n&1&1\\
1&1&n&1&1&n&1&1&1&1&n&1\\
1&1&1&1&1&1&1&1&1&1&1&1\\
\end{array}\right)$$

By computing $Tr(G_{6n}^{-1}G_{6s})$ we get the formula in the
statement.
\end{proof}

The $l=3$ computation was done with Maple. The question of
inverting general Gram matrices is for the moment beyond our
understanding.

In fact, the only conceptual result here is Di Francesco's formula
in \cite{df}, for the determinant of the matrix (which appears as
denominator when computing inverses) in the simplest case, namely
when we have all Temperley-Lieb diagrams.

\section{Asymptotic freeness}

The results in previous section show that the whole subject
belongs to analysis: some $0,1,2,3,\ldots$ hierarchy of problems
is required, in order to make further advances.

In the general context of integration problems, this kind of
revelation first appeared to Weingarten \cite{we}. His idea was to
regard the various formulae as power series in $n^{-1}$, the
inverse of the dimension, and to look at small order coefficients.
In our specific quantum group situation, it is known from
\cite{bc1}, \cite{bc2} that the reasonable quantity to look at is
the asymptotic law of $u_{11}+\ldots+u_{ss}$, with $s\simeq tn$.

\begin{proposition}
We have the asymptotic moment formula
$$\lim_{n\to\infty}\int (u_{11}+\ldots +u_{ss})^{k}=
\sum_p t^{|p|}$$
with $s=[tn]$, where the sum is over all partitions in $NC_h(k)$.
\end{proposition}

\begin{proof}
Consider the diagonal matrix $\Delta_{kn}$ formed by the  diagonal
entries of $G_{kn}$. We have the following formula:
\begin{eqnarray*}
(\Delta_{kn}^{-1/2}G_{kn}\Delta_{kn}^{-1/2})(p,q)
&=&\Delta_{kn}^{-1/2}(p,p)G_{kn}(p,q)\Delta_{kn}^{-1/2}(q,q)
\cr&=&n^{-|p|/2}n^{|p\vee q|}n^{-|q|/2}\cr&=&n^{|p\vee
q|-(|p|+|q|)/2}
\end{eqnarray*}

The exponent is negative for $p\neq q$, and zero for $p=q$. This
gives:
\begin{eqnarray*}
G_{kn}&=&\Delta_{kn}^{1/2}(Id+O(n^{-1/2})) \Delta_{kn}^{1/2}\cr
G_{kn}^{-1}&=&
\Delta_{kn}^{-1/2}(Id+O(n^{-1/2}))\Delta_{kn}^{-1/2}
\end{eqnarray*}

We obtain the following moment estimate:
\begin{eqnarray*}
\int (u_{11}+\ldots +u_{ss})^{k} &=&Tr(G_{kn}^{-1}G_{ks})\cr
&\simeq&Tr(\Delta_{kn}^{-1}\Delta_{ks})\cr &=& \sum_p
n^{-|p|}s^{|p|}
\end{eqnarray*}

With $s=[tn]$ we get the formula in the statement.
\end{proof}

\begin{theorem}
The asymptotic measures
$$\mu_t=\lim_{n\to\infty}{\rm law}(u_{11}+\ldots +u_{[tn][tn]})$$
form a one-parameter truncated semigroup with respect to free
convolution.
\end{theorem}

\begin{proof}
The idea is to follow the proof from the classical case, namely for the group $H_n$, with the log of the Fourier transform replaced by Voiculescu's $R$-transform.

We denote by $F_{kb}$ of number of partitions in $NC_h(2k)$ having $b$ blocks. We set $F_{kb}=0$ for other integer values of $k,b$. All sums will be over integer indices $\geq 0$.

With these notations, the above result shows that the generating series of the moments of $\mu_t$ is given by $f(z)=g(z^2)$, where:
$$g=\sum_{kb}F_{kb}z^kt^b$$

For a partition in $NC_h(2k+2)$, consider the last two legs of the first block. These correspond to certain numbers in $\{1,2,\ldots,2k+1\}$, that we denote $2x_0+1$ and $2(x_0+x_1)+2$. With the notation $x_2=k-x_0-x_1$, we see that our partition is obtained from a partition in $NC_h(2x_0)$, a partition in $NC_h(2x_1)$, and a partition in $NC_h(2x_2)$.

Thus the numbers $F_k=\# NC_h(2k)$ satisfy the following relation:
$$F_{k+1}=\sum_{\Sigma x_i=k} F_{x_0}F_{x_1}F_{x_2}$$

In terms of the numbers $F_{kb}$, this can be written in the following way:
$$\sum_bF_{k+1,b}=\sum_{\Sigma x_i=k}\sum_{u_i}F_{x_0u_0}F_{x_1u_1}F_{x_2u_2}$$

In this formula, each $F_{x_0u_0}F_{x_1u_1}F_{x_2u_2}$ term contributes to $F_{k+1,b}$ with $b=\Sigma u_i$, except for $F_{00}F_{x_1u_1}F_{x_2u_2}=F_{x_1u_1}F_{x_2u_2}$, which contributes to $F_{k+1,b+1}$. We get:
\begin{eqnarray*}
F_{k+1,b}&=&\sum_{\Sigma x_i=k}\sum_{\Sigma u_i=b}F_{x_0u_0}F_{x_1u_1}F_{x_2u_2}\cr
&+&\sum_{\Sigma x_i=k}\sum_{\Sigma u_i=b-1}F_{x_1u_1}F_{x_2u_2}\cr
&-&\sum_{\Sigma x_i=k}\sum_{\Sigma u_i=b}F_{x_1u_1}F_{x_2u_2}
\end{eqnarray*}

By multiplying by $t^b$ and summing over $b$, this gives the following formula for the polynomials $P_k=\sum_bF_{kb}t^b$:
$$P_{k+1}=\sum_{\Sigma x_i=k}P_{x_0}P_{x_1} P_{x_2}+(t-1)\sum_{\Sigma x_i=k}P_{x_1}P_{x_2}$$

In terms of $g=\sum_kP_kz^k$, we have the following equation:
$$g-1=zg^3+(t-1)zg^2$$

In terms of the series $f(z)=g(z^2)$, we have:
$$f-1=z^2f^3+(t-1)z^2f^2$$

In other words, $f$ satisfies the following equation:
$$f=1+(zf)^2(f+t-1)$$

We can turn now to the last part of the proof, which uses Voiculescu's $R$-transform, and general theory from free probability \cite{vdn}. The $R$-transform of a real probability measure $\mu$ is constructed as follows:
\begin{enumerate}
\item $f(z)$ is the generating series of the moments of $\mu$.
\item $G(\xi)=\xi^{-1}f(\xi^{-1})$ is the Cauchy transform.
\item $K(z)$ is defined by $G(K(z))=z$.
\item $R(z)=K(z)-z^{-1}$.
\end{enumerate}

In terms of the Cauchy transform, the equation of $f$ translates as follows:
$$\xi G=1+\xi G^3+(t-1)G^2$$

Thus the $K$-transform satisfies the following equation:
$$zK=1+z^3K+(t-1)z^2$$

This gives the following equation for the $R$-transform:
$$1+zR=1+tz^2+z^3R$$

We can solve this equation, with the following solution:
$$R=\frac{tz}{1-z^2}$$

Thus the $R$-transform is linear in $t$, and this gives the
result.
\end{proof}

Observe that the above proof also reads that the $k$-th free cumulant of $\mu_t$ is $t$ if $k$ is even, and
zero if $k$ is odd. This fact can also be seen directly from combinatorial considerations together with the free
moment-cumulant formula. So, it is also possible to write a proof of Theorem 10.2 of combinatorial flavour.

As a last comment, the set $NC_h(2k)$ is the $s=2$ particular case of the set $NC^{(s)}(sk)$, consisting of partitions of $\{1,2,\ldots,sk\}$ into blocks of size multiple of $s$. These sets, introduced by Edelman \cite{ed}, have been studied in connection with several situations, and have many interesting properties. See e.g. Armstrong \cite{ar} and Stanley \cite{st}.

The general results about $NC^{(s)}(sk)$ can be used at $s=2$ in order to get further information about the semigroup $\{\mu_t\}$. As an example, the numbers $F_{kb}$ in the proof of the above theorem are known to be the Fuss-Narayana numbers, namely:
$$F_{kb}=\frac{1}{b}\begin{pmatrix}k-1\cr b-1\end{pmatrix}\begin{pmatrix}2k\cr b-1\end{pmatrix}$$

This gives the moment formula announced in the introduction.

It is not clear whether at $s\geq 3$ the combinatorics of $NC^{(s)}(sk)$ corresponds or not to some quantum group. The technical problem is actually at the level of probability measures, where positivity seems to be lacking at $s\neq 1,2$. However, at a purely combinatorial level there is no problem, and a one-parameter semigroup of ``free Bessel laws'' seems to exist for any $s$, at least in some virtual sense. These questions are to be added to those raised in the end of the present paper. 

\section{Free quantum groups: the orthogonal case}

We present here some generalizations of notions and results in
previous sections. This will lead to a kind of
axiomatization for free quantum groups, improving previous
considerations from \cite{bc2}. The starting point is the
following conclusion.

\begin{theorem}
Let $A$ be one of the series $A_s(n),A_h(n),A_o(n)$. Let $A_{com}$
 be the series of maximal commutative quotients, and consider the
following measures:
$$\mu_t=\lim_{n\to\infty}{\rm law}(u_{11}+\ldots +u_{[tn][tn]})$$
\begin{enumerate}
\item For $A_{com}$ these form a truncated semigroup with respect
to convolution.\item For $A$ these form a truncated
semigroup with respect to free convolution.
\end{enumerate}
\end{theorem}

This statement was pointed out in \cite{bc2} in the symmetric and
orthogonal case, with the comment that this is a 
formulation of ``freeness''. In other words, what we have shown in
this paper is that the series $A_h(n)$ is free in the sense of
\cite{bc2}. Actually, the story of this paper is that we wanted to
construct a new example, and by starting with the above freeness
notion we eventually ended up with the cube and with $A_h(n)$.

One thing to be done now is to look back at the proof of the above
result, and try to understand where freeness comes from, and axiomatize it. A first notion that emerges from this kind of study
is as follows.

\begin{definition}
A Hopf algebra realizing a factorization
$$A_o(n)\to A\to A_s(n)$$
of the canonical map from left to right is called homogeneous.
\end{definition}

As a first remark, all three algebras $A_o(n),A_h(n),A_s(n)$ are
homogeneous.

Consider now the following Tannakian duality diagram, where $C$
symbols denote as usual categories generated by tensor powers of
the fundamental corepresentation:
$$\begin{matrix}
A_o(n)&\to&A&\to&A_s(n)\cr\cr
\updownarrow&&\updownarrow&&\updownarrow\cr\cr
CA_o(n)&\subset&CA&\subset&CA_s(n)
\end{matrix}$$

In the lower right corner we have the category spanned by all
Temperley-Lieb diagrams. Moreover, in the studied cases $A=A_x(n)$
with $x=o,h,s$ we know that $CA$ is spanned by certain
Temperley-Lieb diagrams. This suggests the following definition.

\begin{definition}
An homogeneous Hopf algebra is called free if the corresponding
tensor category is spanned by Temperley-Lieb diagrams.
\end{definition}

In other words, we ask for the existence of a set of diagrams
$TL\subset TL_s$ such that the above Tannakian duality diagram
becomes:
$$\begin{matrix}
A_o(n)&\to&A&\to&A_s(n)\cr\cr
\updownarrow&&\updownarrow&&\updownarrow\cr\cr
TL_o(n)&\subset&TL&\subset&TL_s(n)
\end{matrix}$$

We denote by $NC$ the set of non-crossing partitions associated to
arches in $TL$. The general techniques in previous sections apply,
and give the following results.
\begin{enumerate}
\item The Weingarten formula holds, with $G_{kn}$ being the Gram
matrix of $NC$. \item The trace formula for sums of diagonal
coefficients holds as well.\item We don't have the asymptotic
moment formula, nor asymptotic freeness.
\end{enumerate}

The problem with asymptotic results is that we don't have any
asymptotics at all. This would require series of free Hopf
algebras $A_x(n)$ having the property that the corresponding sets
of partitions are equal to a given set $NC_x$, for $n$ big enough.
For the moment, we don't have further examples in this sense.

\section{Free quantum groups: the unitary case}

The considerations in previous section can be further generalized. The idea is that
the various notions extend to the unitary case,
according to the following rules:
\begin{enumerate}
\item The orthogonal matrices are replaced by biunitary matrices.
These are square matrices satisfying $u^*=u^{-1}$ and
$\bar{u}=(u^t)^{-1}$. See Wang \cite{wa1}.
 \item The orthogonal Hopf algebras are replaced by unitary Hopf
 algebras. The fundamental corepresentation $u$ must be biunitary, the formulae for $\Delta,\varepsilon$ are
 the same, and the antipode is given by $S(u_{ij})=u_{ji}^*$. See Woronowicz \cite{wo1}.
 \item The algebra $A_o(n)$ is replaced by the algebra $A_u(n)$,
 generated by entries of a $n\times n$ biunitary matrix. A Hopf
 algebra realizing a factorization $A_u(n)\to A\to A_s(n)$ of the
 canonical map from left to right is
 called homogeneous.
 \item The tensor categories generated by $u$ are replaced by
 tensor categories generated by $u,\bar{u}$. The freeness
 condition is that $CA\subset CA_s(n)=TL_s(n)$ is spanned by
 Temperley-Lieb diagrams. We use here the new meaning of $TL_s(n)$.
\end{enumerate}

In other words, we have the following definition in the unitary case.

\begin{definition}
A free Hopf algebra is a quotient of $A_u(n)$, having $A_s(n)$ as quotient, and having the property that its tensor category is spanned by Temperley-Lieb diagrams.
\end{definition}

The fact that we have a Weingarten formula and a moment formula in
this general situation follows from \cite{bc1}, where the case of
$A_u(n)$ is investigated. That is, the results follow by combining
the arguments in previous sections with those in \cite{bc1}.

As in the orthogonal case, the construction of examples is a quite
delicate problem. However, we have here the following simple construction. Given an orthogonal Hopf algebra $(A,u)$, we can define a 
unitary Hopf algebra $(\tilde{A},\tilde{u})$ in the following way:
\begin{enumerate}
\item $\tilde{A}\subset C^*(Z)*A$ is the subalgebra generated by
entries of $\tilde{u}$.\item $\tilde{u}=zu$, where $z$ is the
canonical generator of $C^*(Z)$.
\end{enumerate}

This notion is introduced in \cite{ba}, with the remark that we have
$A_u(n)=\tilde{A}_o(n)$. It follows from results in there that if
$A$ is free orthogonal, then $\tilde{A}$ is free unitary. Thus we
have the following examples of free Hopf algebras:
$$\begin{matrix}
A_u(n)&\to&\tilde{A}_h(n)&\to&\tilde{A}_s(n)\cr\cr
\downarrow&&\downarrow&&\downarrow\cr\cr
A_o(n)&\to&A_h(n)&\to&A_s(n)
\end{matrix}$$

We don't know what is the precise structure of the upper right
algebras. We don't have either free unitary examples not coming
from free orthogonal ones.

Summarizing, we have very few examples of free Hopf algebras.

This suggests that free Hopf algebras might be classifiable by discrete data, say analogous to Coxeter-Dynkin diagrams classifying simple Lie
groups. We don't know if it is so, but we intend to come back to this question in future work.

We should mention that such a project might look a bit too ambitious. It is most likely that a complete classification of free Hopf algebras, if any, would rather have the size of a book. Probably most illustrating here is the story of the table in the introduction: each column is based on a certain number of journal articles, and complexity grows with each column added. So, the hope would be that the list is not too long.

As a conclusion, the reasonable question that we would like to raise here is that of finding a 4-th free quantum group, in the orthogonal case. The answer would probably require some new ideas, even at level of various combinatorial ingredients.

\end{document}